\documentclass[12pt]{article}
\usepackage{graphicx}
\usepackage{epstopdf} 
\usepackage{amssymb} 
\usepackage{amsmath} 

\DeclareMathOperator{\Tr}{Tr}

\numberwithin{equation}{section}
\newtheorem{thm}{Theorem}[section]

\begin{document}

 \centerline{\LARGE Counting covering cycles \vspace{3mm} }
 \centerline{\large G.A.T.F.da Costa 
\footnote{g.costa@ufsc.br} and M. Policarpo}
\centerline{\large Departamento
de Matem\'{a}tica} 
\centerline{\large Universidade Federal de
Santa Catarina} \centerline{\large
88040-900-Florian\'{o}polis-SC-Brasil}

\begin{abstract}

We compute the number of equivalence classes of nonperiodic covering cycles of given length in a non oriented connected graph. A covering cycle is a closed path that traverses
each edge of the graph at least once. A special case  is the number of Euler cycles in the non oriented graph. 
An identity relating the numbers of covering cycles of any length in a graph to a product of determinants is obtained.

\vspace{5mm}

\end{abstract}

\section{ Introduction}

Let $G=(V,E)$ be a finite  connected non oriented graph, $V$ is the set of
$|V|$ vertices, $E$ is the set of edges with $|E|$ elements  labelled $e_{1}$,
...,$e_{|E|}$. The
graph may have multiple edges and loops.
Let's define a new graph from $G$ as follows. Firstly, orient the edges. The orientation is arbitrary but fixed. Let $G'$ be the graph with $2|E|$ oriented
edges built from $G$ by adding
in the opposing oriented edges $e_{|E|+1}=(e_{1})^{-1}$,
...,$e_{2 |E|}=(e_{|E|})^{-1}$, $(e_{i})^{-1}$ being the oriented edge
opposite to $e_{i}$ and with origin (end) the end (origin) of
$e_{i}$. In the case that $e_{i}$ is an oriented loop,
$e_{i+|E|}=(e_{i})^{-1}$ is just an additional oriented loop hooked
to the same vertex.
A path in $G$ is an ordered sequence of edges  in  $G'$,  $(e_{i_{1}},...,e_{i_{N}})$,
$e_{i_{k}} \in E'=\{e_{1}, ..., e_{2E} \}$,
 such that the end of $e_{i_{k}}$ is the origin of $e_{i_{k+1}}$. 
In this paper, all paths are {\it cycles}. A cycle is a non-backtracking tail-less closed path, that is, the end of $e_{i_{N}}$ coincides with the origin of $e_{i_{1}}$, subjected to the non-backtracking condition that $e_{i_{k+1}}  \neq e_{i_{k}+|E|}$. In another words,
a cycle  never goes immediately backwards over a  previous edge. Tail-less means that $e_{i_{1}}  \neq e_{i_{N}}^{-1}$.
The length of a cycle is the number of edges in its sequence. A cycle $p$ is called periodic if $p=q^r$ for some $r>1$, the integer
$r$ is the period of $p$, and
$q$ is a non periodic cycle. The cycle $(e_{i_{N}}, e_{i_{1}}, ...,e_{i_{N-1}})$ is called a circular permutation of $(e_{i_{1}},...,e_{i_{N}})$. Cycles that are circular permutations of one another are said to be equivalent. The cycle
 $(e_{i_{N}}^{-1},...,e_{i_{1}}^{-1})$ is an inversion of $(e_{i_{1}},...,e_{i_{N}})$.
A cycle and its inverse are taken as distinct.

The number of cycles of given length in a non oriented graph $G$ can be computed using the {\it the edge adjacency matrix} of $G$. This 
is the $2|E| \times 2|E|$ matrix $T$ defined as
$T_{ij}=1$, if end vertex of edge $i$ is the origin of
edge $j$ and edge $j$ is not the inverse edge of $i$;
$T_{ij}=0$, otherwise. The number $Tr T^{N}$ counts  backtrack-less and tail-less cycles
of length $N$ in $G$.  The edge adjacency matrix is usually associated with {\it the Ihara zeta function} of a non oriented graph [9].

 {\it A cycle is called an edge covering cycle if  every edge in the graph $G$ is present at least once in the cycle}. Therefore, if $N$ is the length of an edge covering cycle, $N \geq |E|$.
 We are interested in counting equivalence classes of
edge covering cycles in a non oriented graph. For that reason we consider graphs not having degree one vertices although proper subgraphs may have them.  Restricted to this subgraph, say $g$, $Tr T_{g}^{N}$ will count only the backtrack-less and tail-less cycles of length $N$ in $g$ so that edges with 1-degree vertices may be deleted.

In section 2, we obtain a relation for counting the number of equivalence classes of non periodic covering cycles of given length in a non oriented graph.  A special case is the number of Euler cycles in the
non oriented graph. This was our initial motivation. But then we decided to consider  covering cycles in general in order to find the determinantal identity in Theorem 2.4, inspired by an analogous but simpler identity associated to the Ihara zeta function of a graph. The identity resembles some identities associated with dynamical zeta functions [8]. The number of Euler cycles can be computed from one of the coefficients in the expanded form of this identity. 
We give examples in section 3.

\section{ Counting covering and Euler cycles}

Some years back I asked Prof. Michael Krebs (State University of California) [6] if he knew how, along the lines of the combinatorial trace method [5], to compute the number of covering cycles of a given length in a graph. Few days later he replied and his ideas are encapsulated in Theorem 2.1.  

\begin{thm}
Denote by $\{G_{k}\}$ the set of graphs obtained from a graph $G$ as the result of deleting $k$ of its  distinct edges, and by $T_{G_{k}}$ the edge adjacency  matrix of $G_{k}$. Set $G_{0}:=G$ and $G_{|E|}$ the empty graph with 
$Tr T_{G_{|E|}}=0$.  The number of covering cycles of length $N$ in $G$ is
\begin{equation}
\omega(N)= \sum_{k=0}^{|E|} (-1)^{k} \sum_{\{{G}_{k}\} } \Tr {T}_{{G}_{k}}^{N}.
\end{equation}
\end{thm}

\noindent {\bf Proof:}
Apply the inclusion-exclusion principle to find the number of cycles that do not traverse edge $e_{1}$, or $e_{2}$, ..., or $e_{|E|}$. Then, subtract the result from $\Tr T^{N}$ to get (2.1).
$\Box$

\noindent{\bf Remark 2.1.}
We have found in the literature some results that resembles very much (2.1) but associated to the counting of hamiltonian paths. They employ the vertex adjacency matrix. See [1], [4] and [7].

\noindent{\bf Remark 2.2.}
If $N < |E|$, $\Tr T^{N}$ counts the non covering cycles of length $N$ in the graph, but this is exactly what the other terms
in the right hand side of (2.1)  are doing in the formula so we get $\omega=0$ in this case, as expected. If $N \geq |E|$,
it can happen  that $\omega(N)=0$. 
For instance, in a bipartite graph there is no covering cycle  of odd length, hence, only the non covering cycles are counted by $Tr T^{N}$  in this case but, again, this is exactly what the other terms
 together are doing in the formula so we get $\omega=0$ in this case, too.

\begin{thm}
Let $\Theta(N)$ be the number of  equivalence classes of  all   non periodic covering cycles of 
 length $N$ in the graph $G$.  Then,
\begin{equation} \label{b}
\Theta(N) = \frac{1}{N}
 \sum_{g|N} \mu(g) \hspace{1mm} \omega \left(\frac{N}{g} \right)
\end{equation}
The number of Euler cycles in $G$ is $\Theta(|E|)/2=\omega(|E|)/2|E|$. 
\end{thm}

\noindent{\bf Proof:}
In the set of $\omega(N)$  covering cycles  there
is the subset with  $N \Theta(N)$ elements formed by the non periodic covering cycles of length $N$ plus their circular
permutations 
and the subset with  
\begin{equation*}
\sum_{g \neq 1|N}\frac{N}{g} \Theta \left( \frac{N}{g} \right)
\end{equation*}
elements formed by the periodic covering cycles
 of length $N$
(whose periods are the common divisors of $N$) plus their circular
permutations. 
A covering cycle of period $g$ and length $N$  is of
the form
\begin{equation*}
p^{g}=({e}_{{k}_{1}} {e}_{{k}_{2}}
...{e}_{{k}_{\alpha}})^{g}
\end{equation*}
where $\alpha= N/g$, and
$p$ is a non periodic cycle
so that the number of periodic cycles with period $g$ plus their
circular permutations  is given by 
\begin{equation*}
\frac{N}{g} \Theta \left( \frac{N}{g} \right).
\end{equation*}
Hence,
\begin{equation}
\omega(N)= \sum_{g \mid N} \frac{N}{g} \hspace{1mm} \Theta \left( \frac{N}{g}
\right)
\end{equation}
M\"obius inversion formula gives the result. An Euler cycle is a 
cycle that covers each edge of $G$ only once so the number of Euler cycles is given by $\Theta(|E|)$. Divide this number by $2$ to avoid inversions. Notice that  $\omega(|E|/g)=0$ if $g > 1$ by Remark 2.2 so that $\Theta(|E|)=\omega(|E|)/|E|$. Of course, $\Theta(|E|)$ may be zero. This will be the case if the graph has vertices of odd degree.
$\Box$

\begin{thm}
Define
\begin{equation}
h(z):=\sum_{N=1}^{\infty} \frac{ {\omega}(N) }{N} z^{N}.
\end{equation}
Then,
\begin{align}
\prod_{N=1}^{+\infty}  (1-z^{N})^{ \pm \Theta(N)}
= 
e^{\mp  h(z)}  &=
\prod_{k=0}^{|E|} \prod_{\{{G}_{k}\} } \left[  \det(1-z {T}_{{G}_{k}}) \right]^{\pm (-1)^{k}}\\
&= 1 \mp \sum_{i=1}^{+\infty} {d}_{\pm}(i) z^{i}
\end{align}
and
\begin{equation} \label{h}
{d}_{\pm}(i)= \sum_{m=1}^{i} \lambda_{\pm}(m) \sum_{
\begin{array}{l} {a}_{1}+2{a}_{2}+...+i{a}_{i} =i\\
{a}_{1}+...+{a}_{i} = m \end{array}} 
 \prod_{k=1}^{i} 
\frac{(\omega (k))^{{a}_{k}}}{{a}_{k}! k^{{a}_{k}}}
\end{equation}
with $\lambda_{+}(m)=(-1)^{m+1}$, $\lambda_{-}(m)=+1$.
\end{thm}

\noindent{\bf Proof:} Compute the formal logarithm of the infinite product to prove the first equality in (2.5). Then, use (2.1) to get the second. (2.10) follows from Fa\'a di Bruno relation for the exponential. See [2].
$\Box$

\begin{thm}
Given a graph with $n_{0}:=|E|$ edges,
\begin{align*}
a) &\hspace{2mm}  n {d}_{\pm} (n) =  \omega(n) \mp \sum_{k=1}^{n-1}  \omega(n-k) {d}_{\pm}(k), n \geq n_{0},\\
b)& \hspace{2mm}  {d}_{-}(n) = {d}_{+}(n)+\sum_{i=1}^{n-1} {d}_{+}(i) {d}_{-}(n-i), n \geq n_{0},\\
c) & \hspace{2mm} {d}_{\pm}(n)=0, \hspace{2mm}  n \leq {n}_{0}-1,\\
d) & \hspace{2mm} {d}_{\pm}(n)=\frac{\omega(n)}{n}, \hspace{2mm} {n}_{0} \leq n<2{n}_{0},\\
e) & \hspace{2mm}  {d}_{+}(n) \leq \frac{\omega(n)}{n}, \hspace{2mm} 2{n}_{0} \leq n < 3{n}_{0},\\
f) & \hspace{2mm} {d}_{-}(n) \geq \frac{\omega(n)}{n}, \hspace{2mm} n \geq 2{n}_{0},\\
g) & \hspace{2mm} |{d}_{+}(n)| \leq {d}_{-}(n).
\end{align*}
The number of Euler cycles is  $d_{+}(|E|)/2$.
\end{thm}

\noindent{\bf Proof:} See [3] for the proof of a) and b). The other items are an easy consequence of (2.7), a) and b) and Remark 2.2.

\noindent{\bf Remark 2.3.}
If $d_{+}(n) \geq 0$ for all $n$ then  a) implies that $d_{+}(n) \leq \omega(n)/n$ for all $n$. However, positivity does not hold for all graphs. For instance, the cycle graph which has one face and $|E|=|V|=n$ edges. For this graph, the double product of determinants in (2.5) is reduced to
\begin{equation}
\det (1-z T) = (1-z^{|V|})^{2}= 1-2 z^{|V|}+z^{2 |V|}
\end{equation}
which is the reciprocal of the Ihara zeta function of the graph.
This seems to be the only case but a proof is lacking at the moment.

\section{ Examples}

\noindent{\bf Example 1.}
$G_{1}$, the non oriented graph with $R \geq 2$ edges  hooked to a single vertex. Let's orient the edges counterclockwisely. The edge matrix for $G_{1}$ is the $2R \times 2R$ symmetric matrix defined on $G_{1}'$, 
\begin{equation*}
{T}_{{G}_{1}} = \left( \begin{array}{clcr}
A & B\\
B & A 
\end{array} \right)
\end{equation*}
where $A$ is the $R \times R$ matrix with all entries equal to $1$ and $B$  is the $R \times R$ matrix with  the main diagonal entries  equal to $0$ and all the other entries equal to $1$. This matrix has trace given by
\begin{equation*}
\Tr {T}_{{G}_{1}}^{N} = 1+(R-1)(1+(-1)^{N})+(2 R-1)^{N}, \hspace{2mm} N=1,2, \dots
\end{equation*}
Denote by $G_{k}'$ the graph obtained from $G_{1}$ deleting $k$ edges. Then,
\begin{equation*}
\Tr {T}_{{G}_{k}'}^{N} = 1+(R-k-1)(1+(-1)^{N})+(2 R-2 k-1)^{N}, \hspace{2mm} N=1,2, \dots
\end{equation*}
and
\begin{align*}
\omega_{{G}_{1}}(N)&=
\sum_{k=0}^{R-1} (-1)^{k} \left( \begin{array}{cl}
R \\
k  
\end{array} \right)[1+(R-k-1)(1+(-1)^{N})+(2R-2k-1)^{N}] 
\end{align*}
In the special case $R=2$,
the number of classes of nonperiodic covering cycles of length $N$ is
\begin{equation*} 
\Theta_{{G}_{1}}(N) = \frac{1}{N}
 \sum_{g|N} \mu(g) \hspace{1mm} \left( (-1)^{\frac{N}{g}}-2+3^{\frac{N}{g}} \right)
\end{equation*}
$\Theta_{G_{1}}(2)/2=2$ is the number of Euler cycles.
Furthermore,
\begin{align*} \label{g}
\prod_{N=1}^{+\infty}  (1-z^{N})^{\Theta(N)}
&=
\prod_{k=0}^{2} \prod_{\{{G}_{k}\} } \left[  \det(1-z {T}_{{G}_{k}}) \right]^{ (-1)^{k}}\\
&=\frac{(1-2 z- 3 z^2)  }{(1-z)^{2}    }= 1-\sum_{n=1}^{+\infty} 4(n-1)z^n
\end{align*}

\noindent{\bf Example 2.}
$G_{2}$, the bipartite graph with two vertices and three edges with the same orientation linking them. 
The edge matrix $T_{G_{2}}$, defined on $G_{2}'$, is
\begin{equation*}
{T}_{{G}_{2}} = \left( \begin{array}{clcrclcr}
0 & 0 & 0 & 0 & 1 & 1\\
0 & 0 & 0 & 1 & 0 & 1 \\
0 & 0 & 0 & 1 & 1 & 0 \\
0 & 1 & 1 & 0 & 0 & 0\\
1 & 0 & 1 & 0 & 0 & 0\\
1 & 1 & 0 & 0 & 0 & 0 
\end{array} \right) 
\end{equation*}
The matrix has the trace given by
$Tr T_{G_{2}}^{N}= 0$ if $N$ is odd and $Tr T_{G_{2}}^{N}=
4+2 \cdot 2^{N}$ if $N$ is even, so
$\omega_{G_{2}}(N)=0$, if $N$ is odd, and 
\begin{equation*}
\omega_{{G}_{2}}(N)=2 \cdot 2^{N} - 8,
\end{equation*}
 if $N$ is even. 
In this case the number of nonperiodic cycles of length N is
$\Theta_{G_{2}}(N)=0$, if $N$ is odd, hence, there are no Euler cycles, and
\begin{equation*} 
\Theta_{G_{2}}(N) = \frac{1}{N}
 \sum_{g|N,N/g \hspace{1mm}even} \mu(g) \hspace{1mm} \left[ 2 \cdot 2^{\frac{N}{g}}-8 \right],
\end{equation*}
We have
\begin{align*}
\prod_{N=1}^{+\infty} (1-z^{N})^{\Theta_{{G}_{2}}(N)}
&= \frac{(1-4z^{2})}{(z^{2}-1)^{4}}\\
&= 1-\frac{1}{2}\sum_{N=2}^{+\infty}(N+2)(N+1)(N-1) z^{2N}\\
\end{align*}

\section{Final comments}

Let $G=(V,E)$ be a finite  connected and directed graph with no 1-degree vertices,
$|V|$ vertices, and $|E|$ edges $e_{1}$,
...,$e_{|E|}$. 
The
graph may have multiple edges and loops. 
A path in this case is given by an ordered sequence of oriented edges  $(e_{i_{1}},...,e_{i_{N}})$, $i_{k} \in \{1, ..., |E|\}$, the end of $e_{i_{k}}$ is the origin of $e_{i_{k+1}}$. Paths are always backtrack-less, tail-less and have no inverse. 
The number of cycles of a given length in  $G$ is given by the {\it directed vertex adjacency matrix} $A_{d}(G)$ [ 3].
Label the vertices of $G$, $1$, $2$, ..., $|V|$. Then, $(A_{d})_{ij}$ is the number of edges directed from vertex $i$ to vertex $j$. The number of cycles of length $N$ in $G$ is given by $\Tr A_{d}^{N}$. 
Also, one can use the {\it directed edge adjacency matrix} of $G$ to count cycles of a given length in the directed graph. This 
is the $|E| \times |E|$ matrix $S$ defined as
$S_{ij}=1$, if end vertex of edge $i$ is the start vertex of
edge $j$;
$S_{ij}=0$, otherwise. 
All the results in section 2 for non directed graphs hold in the directed case if one replaces the matrix $T$ by  $S$ or $A_{d}$, $Tr S^{N}=Tr A_{d}^{N}$ [3]. In this case $\Theta(|V|)/2=\omega(|V|)/2|V|$ is the number of hamiltonian cycles and we get the result in [4], section 3.6. Also, see [1] and [7].

\subsection*{Acknowledgments}

G.A.T.F.C is deeply grateful to Professor Michael Krebs (State University of California)
for his generous suggestion of the relation (2.1).  Thanks to Prof.
Asteroide Santana (UFSC) for help with  latex commands and determinants.

\end{document}